\documentclass[11pt]{amsart}    

\usepackage{graphicx}
\usepackage{color}
 
\usepackage{amsfonts,amsmath,latexsym,amssymb,amsthm}
\usepackage{hyperref}
\usepackage{geometry}
\usepackage{txfonts}
\usepackage{color}
\definecolor{Myred}{cmyk}{0.0,1.0,1.0,0.00}
\definecolor{Mypurple}{rgb}{0.5,0.0,0.5}
\newtheorem{theorem}{Theorem}

\newtheorem{remark}{Remark}

\newtheorem{lemma}{Lemma}
 
\begin{document}

\title{Three-dimensional magnetic Schr\"{o}dinger operator with the potential supported in a tube}

\author{Diana Barseghyan}

\address{ Department of Mathematics, Faculty of Science, University of Ostrava\\ 30.~dubna 22, 70103 Ostrava, Czech Republic}
\email{diana.schneiderova@osu.cz}

\author{Juan Bory-Reyes}
  
\address{ESIME-Zacatenco, Instituto Polit\'ecnico Nacional, M\'exico,\\ CDMX. 07738. M\'exico.} 
\email{juanboryreyes@yahoo.com}

\author{Baruch Schneider}
 
\address{ Department of Mathematics, Faculty of Science, University of Ostrava \\ 30.~dubna 22, 70103 Ostrava, Czech Republic}
\email{baruch.schneider@osu.cz}

%\author{Yifan Zhang}
%\address{Department of Mathematics, Faculty of Science, University of Ostrava\\ 30.~dubna 22, 70103 Ostrava, Czech Republic}
%\email{yifan.zhang.s01@osu.cz}


\maketitle

\begin{abstract}
In this paper, we study the following magnetic Schr\"odinger operator in $\mathbb{R}^3$:
\[
H=(i \nabla +A)^2- \tilde{V},
\]
where $\tilde{V}$ is non-negative potential supported over the tube built along a curve which is a local deformation of a straight one, and $B:=\mathrm{rot}(A)$ is a non-zero and local (i.e., a compact supported) magnetic field. 
Based on some new strategies, we first prove that the magnetic field does not change the essential spectrum of this system.
Finally, in the last section of this paper, we establish the sufficient condition such that the discrete spectrum is empty.
\end{abstract}

\section{Introduction}
\label{s:intro}

%%%%%%
Let us consider an infinite curve $\Gamma\subset\mathbb{R}^3$ which has no self-intersections. This means, we have the graph of a function $\Gamma : \mathbb{R} \to \mathbb{R}^3$. It is well known that $\Gamma$ can be parametrized by its arc length, 
$|\dot\Gamma(s)|=1$, and $\Gamma$ is $C^2$-smooth. 
 % ------------- %
This allows us to define a global system of curvilinear coordinates in the vicinity of $\Gamma$. As long as the second derivative of the curve, denoted by $\ddot\Gamma(s)$, is non-zero, the curve has a unique Frenet triad frame, denoted by $(t, n, b)$.
This is no longer true if the second derivative of the curve vanishes. If $\Gamma$ is a straight line, one can choose a constant (inertial) Frenet frame instead. When a global Frenet frame exists in any of the above situations, the curve is said to possess an appropriate Frenet frame. Then, via the Serret-Frenet formulae, the Frenet frame evolves along the curve \cite{K78}
$$ 
\left(\begin{array}{ccc} \dot t \\ \dot n \\ \dot b \end{array}\right)
= \begin{pmatrix} 0 & \gamma & 0 \\ -\gamma & 0 & \tau \\ 0 & -\tau & 0 \end{pmatrix}
\left(\begin{array}{ccc} t \\ n \\ b \end{array}\right),
 $$ %\end{equation}
% ------------- %
where $\gamma,\,\tau$ are the curvature and torsion of $\Gamma$, respectively.

The second ingredient is the support of the potential cross section.  We assume that it is a bounded open 
connected set, $\omega\subset\mathbb{R}^2$, containing the origin of the coordinates. In case when a global Frenet frame exists given a function $\alpha:\mathbb{R}\to\mathbb{R}$, we define the tubular region
 % ------------- %
 \begin{equation}\label{tube}
 \Omega^\alpha_{\Gamma, \omega}:= x(\mathbb{R}\times \omega),\quad x(s, r, \theta):= \Gamma(s) - r\big[ n(s)\cos(\theta-\alpha) + b(s)\sin(\theta-\alpha) \big],
 \end{equation}

Our specific objectives include assuming that the curve $\Gamma$ coincides with the coordinate axis 
$\mathbb{Z}=\{(0, 0, x_3)\}_{x_3\in\mathbb{R}}$ outside the interval $[-c, c]$ for some $c > 0$. Inside the interval $(-c, c)$ curve $\Gamma$  possesses an appropriate Frenet frame. We construct the tubular region $\Omega$ as follows: for $s\in [-c, c]$ let it coincide with $\Omega^\alpha_{\Gamma, \omega}\cap\{|s|\le c\}:= x([-c, c]\times \omega)$ as given by definition (\ref{tube}) and for $s, \,|s|\ge c$ let it be a straight tube $\{|s|\ge c\}\times \omega$. 

The classical theory of a magnetic field, even a local one, can significantly affect the behaviour of waveguide systems, in particular the existence of a geometrically induced discrete spectrum. While a particle confined in a fixed-profile tube with a hard-wall boundary can exist in localized states whenever the tube is bent or locally deformed (and asymptotically straight), cf. \cite{EK15} for a comprehensive review of quantum waveguide theory, the presence of a local magnetic field can destroy such a discrete spectrum. This is a consequence of a Hardy-type inequality proved by Ekholm and Kova\v{r}\'{\i}k \cite{EK05}. In the applied literature, magnetic fields are usually called $B$-fields. 

Another interesting result related to the magnetic Dirichlet Laplacian in thin tubular neighborhoods of curves is presented in \cite{KR14}, where the authors study a norm resolvent convergence in the limit when the cross section of the tube shrinks to a point and the magnetic field is allowed to explode to infinity.

%In this work we consider three-dimensional magnetic Schr\"{o}dinger operators with an attractive potential supported over a tube, which is a local perturbation of a straight tube. 
We are interested in the magnetic Schr\"{o}dinger operators with a compactly supported magnetic field and an attractive potential, 
$\tilde{V}$, over the tube $\Omega$. 
A similar model with a zero magnetic field was considered in \cite{E22}. If the corresponding tube is not straight, but is straight outside a bounded region, or it is at least asymptotically straight in a suitable sense, then the essential spectrum of the Hamiltonian is stable for potentials with some additional properties. On the other hand, the author derived a condition that ensures the existence of a discrete spectrum when a potential channel is bent and twisted under suitable assumptions within the aforementioned bounded region. In this work, we will show that adding a compactly supported magnetic field can destroy these eigenvalues.

To recall some definitions and set the framework, consider the auxiliary function $V\in L^\infty(\mathbb{R}^2)$, such that 
\begin{equation}\label{assump.}
V\ge 0,\quad \mathrm{supp} V\subset \omega.
\end{equation}

Using the polar coordinates in the plane perpendicular to $\Gamma$ inside the interval $(-c, c)$, we define 
\begin{eqnarray}\label{potential}
\tilde{V}:\: \Omega \to\mathbb{R}_+,\quad \tilde{V}(x(s, r, \theta))=V(r, \theta),  \\\label{Hamiltonian} H =(i \nabla +A)^2- \tilde{V}(x)\,; 
\end{eqnarray}
where $A$ is the vector potential corresponding to the magnetic field $B$. Here $(i \nabla +A)^2$ abbreviates the classical differential operator $(i \partial_1 +A_1)^2 + (i \partial_2 +A_2)^2$.

%----------------%
There are several ways to define the operator properly. The simplest way is to identify it with the
unique self-adjoint operator associated with the closed quadratic form
\begin{equation}\label{q}
q(f)=\int_{\mathbb{R}^3}|i\nabla f+ A f|^2\,d x_1\,d x_2\,d x_3- \int_{\mathbb{R}^3}\tilde{V} |f|^2\,d x_1\,d x_2\,d x_3,\quad f\in \mathcal{H}^1(\mathbb{R}^3)\,,
\end{equation}
defined on the Sobolev space $\mathcal{H}^1(\mathbb{R}^3)$, where for $x\in \mathbb{R}^3$ we used the notation $x=(x_1, x_2, x_3)$.
 
From now on, we will focus our attention on the following assumptions:

\begin{itemize} 
\item
the magnetic field $B = \mathrm{rot}\, A\in L_{\mathrm{loc}}^2(\mathbb{R}^3)$ is compactly supported with $\mathrm{supp} (B) \subset\mathcal{B}(0, 2s_0)$, where $\mathcal{B}(0, 2s_0)$ is the three dimensional ball centered at the origin with radius $2s_0$. Moreover, on $\mathcal{B}\left(0, s_0\right)$ where $\mathcal{B}(0, s_0)$ is the three dimensional ball centered at the origin with radius $s_0$, the magnetic field $B$ is constant and equal to the vector $(B_1^0, B_2^0, B_3^0)$;

\item
outside of the strip $\mathbb{R}^2\times\left\{|x_3|\le s_0/\sqrt{2}\right\}$, the curve $\Gamma$ coincides with the coordinate axis $\mathbf{Z}=\{(0, 0, x_3)\}_{x_3\in\mathbb{R}}$;
\item
unlike the work \cite{E22}, where the author made the additional assumption on the twisting of the tubular region and on torsion of the curve $\Gamma$, we make no additional assumptions except that the curve is $C^2$-smooth.
\end{itemize}

To continue the description, we must introduce a two-dimensional auxiliary operator acting in $L^2(\mathbb{R}^2)$,
$$h_V= -\frac{\partial^2}{\partial x_1^2}- \frac{\partial^2}{\partial x_2^2}- V(x_1, x_2),$$ with domain $\mathcal{H}^2(\mathbb{R}^2)$.
 % ------------- %
 Operator $h_V$ has by (\ref{assump.}) a nonempty and finite discrete spectrum with
 % ------------- %
 \begin{equation} \label{infsph}
e := \inf\sigma_\mathrm{disc}(h_V) = \inf\sigma(h_V)\in \big(-\|V\|_\infty, 0\big).
 \end{equation}
 % ------------- %
We are interested in the spectral properties of the operator $H$. Let us mention, that in work \cite{E22}, it was shown that if $\Omega$ is a straight tube and the magnetic field is absent, then the spectrum of the operator $H$ is purely essential and coincides with $[e, \infty)$.The same work shows that local perturbations of a straight tube do not change its essential spectrum. However, the discrete spectrum of such an operator is nonempty if the potential is large and the tube is narrow enough. We will study the operator $H$ with a non-zero magnetic field. 
We are going to prove that the existence of a compactly supported magnetic field does not change the essential spectrum.  On the other hand, we will establish a sufficient condition for the absence of its discrete spectrum.  
As an application, we give the example of paper \cite{E22}, where the curve deformation produces a discrete spectrum and we show that the local magnetic field destroys the eigenvalues.

\section{Main results}\label{main}

\begin{theorem}\label{th.1}
Assuming the conditions stated in the introduction, the essential spectrum of the operator $H$ coincides with the half-line $\big[e, \infty\big)$.
\end{theorem}
 %----------------%
  
\begin{proof}

First, let us show that the essential spectrum of the operator $H$ covers the half-line $[e, \infty)$
\begin{equation}\label{sigmaess}\sigma_{\mathrm{ess}}(H)\supset [e, \infty).
\end{equation}

We are going to use the Weyl criterion \cite[Thm~VII.12]{RS81}. 
Put $\lambda=e+p^2,\,p\in\mathbb{R}$, and consider the sequence of vectors
%----------------%
$$
\psi_k:\: \psi_k(x_1, x_2, x_3)=\frac{1}{\sqrt{k}}\,f(x_1, x_2)\,\mathrm{e}^{i p x_3}\,\chi\left(\frac{x_3}{k}\right),
$$
 %----------------%
where $f$ is the ground state eigenfunction of $h_V$, $\chi$ is a smooth function with support in the interval $(1, 2)$, and $k\in \mathbb{N}$.

It is easy to check that $\psi_k\in \mathrm{Dom}(H)$; idea is to exploit the fact that there is a part of the plane on which the magnetic field has no influence. We can use fact that the divergence of the magnetic field $B$ must to be equal to zero: $\frac{\partial B_1}{\partial x_1}+ \frac{\partial B_2}{\partial x_2}+ \frac{\partial B_3}{\partial x_3}=0$ due to the fact that $\mathrm{rot}(A)=B$ we are able to choose the Landau gauge for the vector potential $A$, i.e. by setting 
\begin{eqnarray}\nonumber
A=(A_1, A_2, A_3)\\\label{gauge}= \left(\int_{2s_0}^{x_3} \int_0^{x_2} \frac{\partial B_1}{\partial x_1}(x_1, t_2, t_3)\,d t_2\,d t_3+ \int_{2s_0}^{x_3} B_2(x_1, x_2, t_3)\,d t_3, 0, \int_0^{x_2} B_1(x_1, t_2, x_3)\,d t_2\right)
\end{eqnarray}
as the solution for $\mathrm{rot}(A)=B$:
\begin{eqnarray}\nonumber
\frac{\partial A_3}{\partial x_2}- \frac{\partial A_2}{\partial x_3}= B_1
\\\nonumber \frac{\partial A_3}{\partial x_1}- \frac{\partial A_1}{\partial x_3}= -B_2
\\\label{system}\frac{\partial A_2}{\partial x_1}- \frac{\partial A_1}{\partial x_2}= B_3\,.
\end{eqnarray}

From this we get $A(x_1, x_2, x_3)=0$ for $x_3> 2s_0$. Then, for sufficiently large values of $k$
\begin{equation}\label{k}\int_{\mathbb{R}^3}|H \psi_k- \lambda \psi_k|^2\,d x_1\,d x_2\,d x_3= \int_{\mathbb{R}^3}|-\Delta \psi_k- \tilde{V} \psi_k- \lambda \psi_k|^2\,d x_1\,d x_2\,d x_3.\end{equation}

One has
\begin{align*}&\frac{\partial^2\psi_k}{\partial x_1^2}=\frac{1}{\sqrt{k}}\,\frac{\partial^2 f}{\partial x_1^2}(x_1, x_2)\,\mathrm{e}^{i p x_3}\,\chi\left(\frac{x_3}{k}\right)\,,\\&
\frac{\partial^2\psi_k}{\partial x_2^2}=\frac{1}{\sqrt{k}}\,\frac{\partial^2 f}{\partial x_2^2}(x_1, x_2)\,\mathrm{e}^{i p x_3}\,\chi\left(\frac{x_3}{k}\right)\,,\\&
\frac{\partial^2\psi_k}{\partial x_3^2}=\frac{1}{\sqrt{k}}\,\left(-p^2 f(x_1, x_2)\,\chi\left(\frac{x_3}{k}\right)+ \frac{2 i p}{k} f(x_1, x_2)\,\chi'\left(\frac{x_3}{k}\right)+ \frac{1}{k^2}
f(x_1, x_2)\,\chi''\left(\frac{x_3}{k}\right)\right) \,\mathrm{e}^{i p x_3}\,.
\end{align*}

Hence in view of (\ref{k}) and the fact that $\tilde{V}(x_1, x_2, x_3)= V(x_1, x_2)$ for $x_3>  s_0/\sqrt{2}$, we have
\begin{gather*}
\int_{\mathbb{R}^3}|H \psi_k- \lambda \psi_k|^2\,d x_1\,d x_2\,d x_3\\\nonumber= \frac{1}{k} \int_k^{2k}\int_{\mathbb{R}^2}\biggl|-\frac{\partial^2 f}{\partial x_1^2}(x_1, x_2)\,\chi\left(\frac{x_3}{k}\right)- \frac{\partial^2 f}{\partial x_2^2}(x_1, x_2)\,\chi\left(\frac{x_3}{k}\right)+ p^2 f(x_1, y_2)\,\chi\left(\frac{x_3}{k}\right)- \frac{2 i p}{k} f(x_1, x_2)\,\chi'\left(\frac{x_3}{k}\right)
\\\nonumber-\frac{1}{k^2} f(x_1, x_2)\,\chi''\left(\frac{x_3}{k}\right)-
V(x_1, x_2)f(x_1, x_2)\,\chi\left(\frac{x_3}{k}\right)- (e+p^2) f(x_1, x_2)\,\chi\left(\frac{x_3}{k}\right)\biggr|^2\,d x_1\,d x_2\,d x_3\,.\end{gather*}

Since $f$ is the ground state eigenfunction of $h_V$ corresponding to the eigenvalue $e$ the above expression implies
\begin{align}\nonumber&
\int_{\mathbb{R}^3}|H \psi_k- \lambda \psi_k|^2\,dx_1\,dx_2\,d x_3\\\nonumber&= \frac{1}{k} \int_k^{2k}\int_{\mathbb{R}^2}\biggr|\frac{2 i p}{k} f(x_1, x_2)\,\chi'\left(\frac{x_3}{k}\right)+ \frac{1}{k^2} f(x_1, x_2)\,\chi''\left(\frac{x_3}{k}\right)\biggr|^2\,d x_1\,d x_2\,d x_3
&\\\nonumber\le &\frac{8 p^2}{k^3} \int_k^{2k}\int_{\mathbb{R}^2} f(x_1, x_2)^2\,\left(\chi'\left(\frac{x_3}{k}\right)\right)^2\,d x_1\,d x_2\,d x_3+ \frac{2}{k^5}  \int_k^{2k}\int_{\mathbb{R}^2} f(x_1, x_2)^2\,\left(\chi''\left(\frac{x_3}{k}\right)\right)^2\,d x_1\,d x_2\,d x_3&\\
\nonumber&= \frac{8 p^2}{k^2} \int_{\mathbb{R}^2}f(x_1, x_2)^2\,d x\,d y \int_1^2 (\chi'(t))^2\,d t+ \frac{2}{k^4} \int_{\mathbb{R}^2} f(x_1, x_2)^2\,d x_1\,d y_2  \int_1^2(\chi''(t))^2\,d t
=\mathcal{O}\left(\frac{1}{k^2}\right),&\end{align}
%----------------%
and since  $\|\psi_k\|^2$ are independent of $k$, we infer that $\lambda =e+p^2\in \sigma(H)$. Moreover, one can choose a sequence $\{k_n\}_{n=1}^\infty$ such that $k_n\to\infty$ as $n\to\infty$ and the supports of the functions $\psi_{k_n}$ are mutually disjoint which means that $e+ p^2\in \sigma_\mathrm{ess}(H)$ for all $p\in \mathbb{R}$, and, one infers that $[e, \infty) \subset \sigma_{\mathrm{ess}}(H)$.

Next, we must show that if the spectrum of $H$ below $e$ exists, it must be discrete. The Neumann bracketing provides an estimate
 %----------------%
\begin{equation}
\label{bracketing}
H\ge H_1\oplus H_2\oplus H_3,
\end{equation}
where $H_1$ is the Neumann restriction of $H$ on $L^2(\mathbb{R}^2\times \{x_3> 2s_0\})$, i.e., the operator with the Neumann boundary conditions at $x_3= 2s_0$, $H_2$ is the Neumann restriction of $H$ on $L^2(\mathbb{R}^2\times \{x_3< -2s_0\})$, i.e., the operator with the Neumann boundary conditions at $x_3= -2s_0$, and $H_3$ is the complementary Neumann restriction on $L^2(\mathbb{R}^2\times\{|x_3|< 2s_0\})$. We will show that the three operators $H_1,\,H_2$ and $H_3$ have purely discrete spectrum below $e$. Then the same will hold for the direct sum $H_1 \oplus H_2\oplus H_3$ and via the minimax principle and (\ref{bracketing}) also for the operator $H$.

Let us start with $H_1$. Now, just by the fact that $\tilde{V}(x_1, x_2, x_3)= V(x_1, x_2)$ and the gauge (\ref{gauge}) $A(x_1, x_2, x_3)=0$ if $x_3>2s_0$, one can easily verify that
 %----------------%
\begin{align*}
 &\int_{\mathbb{R}^2\times \{x_3>2s_0\}} |i \nabla u + A u|^2\,d x_1\,d x_2\,d x_3-  \int_{\mathbb{R}^2\times \{x_3>2s_0\}} \tilde{V}|u|^2\,d x_1\,d x_2\,d x_3 \\ & \ge
\int_{\{x_3>2s_0\}} \int_{\mathbb{R}^2}\left(\left|\frac{\partial u}{\partial x_1}\right|^2+ \left|\frac{\partial u}{\partial x_2}\right|^2\right)\,d x_1\,d x_2\,d x_3+ \int_{\{x_3>2s_0\}} \int_{\mathbb{R}^2}\left|\frac{\partial u}{\partial x_3}\right|^2\,d x_1\,d x_2\,d x_3\\&-  \int_{\{x_3>2s_0\}} \int_{\mathbb{R}^2} V(x_1, x_2) |u|^2\,d x_1\,d x_2\,d x_3 
\\ &
\ge\int_{\{x_3>2s_0\}} \int_{\mathbb{R}^2}\left(\left|\frac{\partial u}{\partial x_1}\right|^2+ \left|\frac{\partial u}{\partial x_2}\right|^2\right)\,d x_1\,d x_2\,d x_3-  \int_{\{x_3>2s_0\}} \int_{\mathbb{R}^2} V(x_1, x_2) |u|^2\,d x_1\,d x_2\,d x_3.\\ 
\end{align*}

Since the principal eigenvalue of the operator $h_V$ is $e$, the above estimate implies
 %----------------%
\begin{align*}
\int_{\mathbb{R}^2\times \{x_3>2s_0\}} |i \nabla u + A u|^2\,d x_1\,d x_2\,d x_3-  \int_{\mathbb{R}^2\times \{x_3>2s_0\}} \tilde{V}(x_1, x_2, x_3)|u|^2\,d x_1\,d x_2\,d x_3&\\\ge e \int_{\mathbb{R}^2\times \{x_3>2s_0\}}|u|^2\,d x_1\,d x_2\,d x_3\,;
\end{align*}
 %----------------%
this means that the spectrum of $H_1$ below $e$ is empty.

The proof that the spectrum of $H_2$ below $e$ is empty can be done in the same way for the operator for operator $H_2$. This time by choosing the solution of system (\ref{system}) as follows:  
\begin{eqnarray}\nonumber
\overline{A}=(\overline{A}_1, \overline{A}_2, \overline{A}_3)=\\ \left(-\int_{x_3}^{-2s_0} \int_0^{x_2} \frac{\partial B_1}{\partial x_1}(x_1, t_2, t_3)\,d t_2\,d t_3- \int_{x_3}^{-2s_0} B_2(x_1, x_2, t_3)\,d t_3, 0, \int_0^{x_2} B_1(x_1, t_2, x_3)\,d t_2\right)\,\nonumber\\
\label{gauge1}
\end{eqnarray}
which satisfies that $\overline{A}(x_1, x_2, x_3)=0$ for $x_3<-2s_0$. 

It remains to deal with the operator $H_3$. Let us choose a positive number, $\alpha$, in such a way that the square, $(-\alpha, \alpha)^2$, contains the set $\Omega\cap\{x_3= z_0\}$ for every $z_0, \,|z_0|< 2s_0$.
Again using the Neumann bracketing method, we get 
\begin{equation}\label{simple}H_3\ge H^1_3 \oplus H_3^2,\end{equation}
where $H_3^1$ is the Neumann restriction of $H_3$ to $L^2((-\alpha, \alpha)^2\times (-2s_0, 2s_0))$ and $H_3^2$ is the Neumann restriction of $H_3$ to $L^2\left(\left(\mathbb{R}^2\setminus (-\alpha, \alpha)^2\right)\times (-2s_0, 2s_0)\right)$.

Let us choose the gauge $A$ as in (\ref{gauge}). Then, the following simple estimate holds: 
$$H_3^1\ge (i \nabla+ A)^2- \|\tilde{V}\|_{L^\infty(\mathbb{R}^3)},$$
which together with the discreteness of the magnetic Neumann Laplacian on the cube $(-\alpha, \alpha)^2\times (-2s_0, 2s_0)$
implies the discreteness of the spectrum of $H_3^1$.

Moving on to the operator $H_3^2$, it is clear that it is non-negative due to the fact that on its domain potential $\tilde{V}$ is zero. Therefore, the negative spectrum of $H_3^2$ is absent.

Therefore, in view of (\ref{simple}), we conclude that the negative spectrum of $H_3$, if it exists, consists of a finite number of eigenvalues of finite multiplicity. Together with (\ref{bracketing}) and the fact that the spectrums of $H_1$ and $H_2$ below $e$ are empty, this completes the proof of our claim.
\end{proof}

\begin{theorem}\label{th.2}
Assume the conditions of Theorem \ref{th.1} hold. Then there exists a constant $C > 0$ such that the discrete spectrum of the operator $H$ is empty, provided that
\begin{eqnarray}
\label{assumption1}
\Omega\cap \left(\mathbb{R}^2\times\{|x_3|\le s_0/\sqrt{2}\}\right) \subset \mathcal{B}(0, s_0)\,,\\\label{assumption2}
\underset{(x_1, x_2),\, (y_1, y_2)\in \omega}{\mathrm{sup}}|V(x_1, x_2) -V(y_1, y_2)|\le C B_3^0\,.
\end{eqnarray}
\end{theorem}

\begin{proof}

We need the following well-known result (see, for example,  \cite[Thm.11.5]{LL01}):
\begin{lemma}\label{Appendix} Let $L:= -\Delta+ Q$, where $Q\le 0$ is a non-zero compactly supported function belonging to $L^\infty(\mathbb{R}^2)$. Then, the operator $L$ has at least one negative eigenvalue. Moreover, the ground state eigenfunction is nowhere zero and can be chosen to be positive.
\end{lemma}

Using the notation from the previous section, let $f$ be the ground state eigenfunction of the operator $h_V$. Given the above lemma, we can choose it to be positive. We can then represent each function, $\psi\in C_0^\infty(\mathbb{R}^3)$ as a multiplication $\varphi(x_1, x_2, x_3) f(x_1, x_2)$ with $\varphi\in C_0^\infty(\mathbb{R}^3)$. Then, with the gauge (\ref{gauge}) for any $\psi\in C_0^\infty(\mathbb{R}^3)$  we have
\begin{eqnarray}\nonumber
\int_{\mathbb{R}^3}\left(|i \nabla \psi(x_1, x_2, x_3)+ A(x_1, x_2, x_3) \psi(x_1, x_2, x_3)|^2- \tilde{V}(x_1, x_2, x_3) |\psi(x_1, x_2, x_3)|^2\right)\,d x_1\,d x_2\,d x_3\\\nonumber=\int_{\mathbb{R}^3}\biggl(\left|i \frac{\partial \varphi}{\partial x_1}(x_1, x_2, x_3) f(x_1, x_2)+  i \varphi(x_1, x_2, x_3) \frac{\partial f}{\partial x_1}(x_1, x_2)+ A_1(x_1, x_2, x_3)\varphi(x_1, x_2, x_3) f(x_1, x_2)\right|^2\\\nonumber+
\left|\frac{\partial \varphi}{\partial x_2}(x_1, x_2, x_3) f(x_1, x_2)+  \varphi(x_1, x_2, x_3) \frac{\partial f}{\partial x_2}(x_1, x_2)\right|^2\\\nonumber-\tilde{V}(x_1, x_2, x_3) |\varphi(x_1, x_2, x_3)|^2 f(x_1, x_2)^2\biggr)\,d x_1\,d x_2\,d x_3\\\nonumber +\int_{\mathbb{R}^3}\left|i \frac{\partial \varphi}{\partial x_3}(x_1, x_2, x_3)+ A_3(x_1, x_2, x_3) \varphi(x_1, x_2, x_3)\right|^2 f(x_1, x_2)^2\,d x_1\,d x_2\,d x_3\\\nonumber = \int_{\mathbb{R}^3}\biggl(\left|\frac{\partial \varphi}{\partial x_1}(x_1, x_2, x_3)\right|^2 f(x_1, x_2)^2+ \frac{\partial \varphi}{\partial x_1}(x_1, x_2, x_3) \overline{\varphi}(x_1, x_2, x_3) f(x_1, x_2) \frac{\partial f}{\partial x_1}(x_1, x_2)\\\nonumber
+ \varphi(x_1, x_2, x_3) \frac{\partial \overline{\varphi}}{\partial x_1}(x_1, x_2, x_3) f(x_1, x_2) \frac{\partial f}{\partial x_1}(x_1, x_2)
\\\nonumber+|\varphi(x_1, x_2, x_3)|^2 \left|\frac{\partial f}{\partial x_1}(x_1, x_2)\right|^2+ i A_1(x_1, x_2, x_3) \frac{\partial \varphi}{\partial x_1}(x_1, x_2, x_3)\overline{\varphi}(x_1, x_2, x_3) f^2(x_1, x_2)\\\nonumber- i A_1(x_1, x_2, x_3)
\frac{\partial \overline{\varphi}}{\partial x_1}(x_1, x_2, x_3) \varphi(x_1, x_2, x_3) f^2(x_1, x_2)+ (A_1)^2(x_1, x_2, x_3)|\varphi(x_1, x_2, x_3)|^2 f^2(x_1, x_2)\\\nonumber
+\left|\frac{\partial \varphi}{\partial x_2}(x_1, x_2, x_3)\right|^2 f(x_1, x_2)^2+ \frac{\partial \varphi}{\partial x_2}(x_1, x_2, x_3) \overline{\varphi}(x_1, x_2, x_3) f(x_1, x_2) \frac{\partial f}{\partial x_2}(x_1, x_2)\\\nonumber
+ \varphi(x_1, x_2, x_3) \frac{\partial \overline{\varphi}}{\partial x_2}(x_1, x_2, x_3) f(x_1, x_2) \frac{\partial f}{\partial x_2}(x_1, x_2)
+|\varphi(x_1, x_2, x_3)|^2 \left|\frac{\partial f}{\partial x_2}(x_1, x_2)\right|^2\\\nonumber
- \tilde{V}(x_1, x_2, x_3) |\varphi(x_1, x_2, x_3)|^2 f(x_1, x_2)^2\biggr)\,d x_1\,d x_2\,d x_3
\\\label{calc.}+\int_{\mathbb{R}^3}\left|i \frac{\partial \varphi}{\partial x_3}(x_1, x_2, x_3)+ A_3(x_1, x_2, x_3)\varphi(x_1, x_2, x_3)\right|^2 f(x_1, x_2)^2\,d x_1\,d x_2\,d x_3\,.
\end{eqnarray}

One can easily verify this by integrating by parts
\begin{eqnarray*}
\int_{\mathbb{R}^3}\biggl(\frac{\partial \varphi}{\partial x_1}(x_1, x_2, x_3) \overline{\varphi}(x_1, x_2, x_3) f(x_1, x_2) \frac{\partial f}{\partial x_1}(x_1, x_2)+ \varphi(x_1, x_2, x_3) \frac{\partial \overline{\varphi}}{\partial x_1}(x_1, x_2, x_3) f(x_1, x_2) \frac{\partial f}{\partial x_1}(x_1, x_2)\\+|\varphi(x_1, x_2, x_3)|^2 \left|\frac{\partial f}{\partial x_1}(x_1, x_2)\right|^2\biggr)\,d x_1\,d x_2\,d x_3= -\int_{\mathbb{R}^3}\frac{\partial^2 f}{\partial x_1^2}(x_1, x_2) f(x_1, x_2) |\varphi(x_1, x_2, x_3)|^2\,d x_1\,d x_2\,d x_3\,,\\
\int_{\mathbb{R}^3}\biggl(\frac{\partial \varphi}{\partial x_2}(x_1, x_2, x_3) \overline{\varphi}(x_1, x_2, x_3) f(x_1, x_2) \frac{\partial f}{\partial x_2}(x_1, x_2)+ \varphi(x_1, x_2, x_3) \frac{\partial \overline{\varphi}}{\partial x_2}(x_1, x_2, x_3) f(x_1, x_2) \frac{\partial f}{\partial x_2}(x_1, x_2)\\+|\varphi(x_1, x_2, x_3)|^2 \left|\frac{\partial f}{\partial x_2}(x_1, x_2)\right|^2\biggr)\,d x_1\,d x_2\,d x_3= -\int_{\mathbb{R}^3}\frac{\partial^2 f}{\partial x_2^2}(x_1, x_2) f(x_1, x_2) |\varphi(x_1, x_2, x_3)|^2\,d x_1\,d x_2\,d x_3\,,
\end{eqnarray*}
which together with (\ref{calc.}) and the fact that $f$ is the ground state eigenfunction of $h_V$ gives
\begin{align}\nonumber&
\int_{\mathbb{R}^3}\left(|i \nabla \psi(x_1, x_2, x_3)+ A(x_1, x_2, x_3) \psi(x_1, x_2, x_3)|^2- \tilde{V}(x_1, x_2, x_3) |\psi(x_1, x_2, x_3)|^2\right)\,d x_1\,d x_2\,d x_3\\\nonumber& = \int_{\mathbb{R}^3}\biggl(\left|\frac{\partial \varphi}{\partial x_1}(x_1, x_2, x_3)\right|^2f(x_1, x_2)^2+ i A_1(x_1, x_2, x_3) \frac{\partial \varphi}{\partial x_1}(x_1, x_2, x_3)\overline{\varphi}(x_1, x_2, x_3) f^2(x_1, x_2)\\\nonumber&- i A_1(x_1, x_2, x_3)
\frac{\partial \overline{\varphi}}{\partial x_1}(x_1, x_2, x_3) \varphi(x_1, x_2, x_3) f^2(x_1, x_2)\\\nonumber&+ (A_1)^2(x_1, x_2, x_3)|\varphi(x_1, x_2, x_3)|^2 f^2(x_1, x_2)\biggr)\,d x_1\,d x_2\,d x_3+
\int_{\mathbb{R}^3}\left|\frac{\partial \varphi}{\partial x_2}(x_1, x_2, x_3)\right|^2 f(x_1, x_2)^2\,d x_1\,d x_2\,d x_3\\\nonumber&+\int_{\mathbb{R}^3}\left|i \frac{\partial \varphi}{\partial x_3}(x_1, x_2, x_3)+ A_3(x_1, x_2, x_3) \varphi(x_1, x_2, x_3)\right|^2 f(x_1, x_2)^2\,d x_1\,d x_2\,d x_3
\\\nonumber&+
\int_{\mathbb{R}^3}\left(\left(-\frac{\partial^2 f}{\partial x_1^2}(x_1, x_2)- \frac{\partial^2 f}{\partial x_2^2}(x_1, x_2) - V(x_1, x_2) f(x_1, x_2)\right)f(x_1, x_2) |\varphi(x_1, x_2, x_3)|^2\right)\,d x_1\,d x_2\,d x_3\\\nonumber&- \int_{\mathbb{R}^3}(\tilde{V}(x_1, x_2, x_3)-V(x_1, x_2)) f(x_1, x_2)^2 |\varphi(x_1, x_2, x_3)|^2\,d x_1\,d x_2\,d x_3
\\\nonumber&
= \int_{\mathbb{R}^3}\left|i \nabla \varphi + A(x_1, x_2, x_3) \varphi(x_1, x_2, x_3)\right|^2 f(x_1, x_2)^2\,d x_1\,d x_2\,d x_3\\\label{est.}&- \int_{\mathbb{R}^3}(\tilde{V}(x_1, x_2, x_3)-V(x_1, x_2)) f(x_1, x_2)^2 |\varphi(x_1, x_2, x_3)|^2\,d x_1\,d x_2\,d x_3+ e \int_{\mathbb{R}^3} |\psi(x_1, x_2, x_3)|^2\,d x_1\,d x_2\,d x_3\,.
\end{align}

Now, let us consider the integral
$$\int_{\mathbb{R}^3}\left|i \nabla \varphi + A(x_1, x_2, x_3) \varphi(x_1, x_2, x_3)\right|^2 f(x_1, x_2)^2\,d x_1\,d x_2\,d x_3.$$ 

Let $$\beta_f:=\underset{\{(x_1, x_2): x_1^2+ x_2^2\le {s_0^2}\}}{\mathrm{min}}\,f(x_1, x_2)^2.$$

From the positivity of $f$, as explained in Lemma \ref{Appendix}, it follows that $\beta_f> 0$. Then

\begin{eqnarray}\nonumber
\int_{\mathbb{R}^3}\left|i \nabla \varphi + A(x_1, x_2, x_3) \varphi(x_1, x_2, x_3)\right|^2 f(x_1, x_2)^2\,d x_1\,d x_2\,d x_3\\\nonumber 
\ge\int_{\{x_1^2+y_2^2+x_3^2\le {s_0^2}\}}\left|i \nabla \varphi(x_1, x_2, x_3) + A(x_1, x_2, x_3) \varphi(x_1, x_2, x_3)\right|^2 f(x_1, x_2)^2\,d x_1\,d x_2\,d x_3
\\\nonumber\ge\int_{\{x_1^2+x_2^2+x_3^2\le {s_0^2}\}}\left| {i\frac{\partial \varphi}{\partial x_1}(x_1, x_2, x_3)+ A_1(x_1, x_2, x_3) \varphi(x_1, x_2, x_3)}\right|^2 f(x_1, x_2)^2\,d x_1\,d x_2\,d x_3
\\\nonumber+\int_{\{x_1^2+x_2^2+x_3^2\le {s_0^2}\}} \left| {\frac{\partial \varphi}{\partial x_2}(x_1, x_2, x_3)}\right|^2 f(x_1, x_2)^2\,d x_1\,d x_2\,d x_3
\\\nonumber\ge\beta_f \int_{\{x_1^2+x_2^2+x_3^2\le {s_0^2}\}}\left|{i \frac{\partial \varphi}{\partial x_1}(x_1, x_2, x_3)+A_1(x_1, x_2, x_3) \varphi(x_1, x_2, x_3)}\right|^2 \,d x_1\,d x_2\,d x_3\\\nonumber
+\beta_f \int_{\{x_1^2+x_2^2+x_3^2\le {s_0^2}\}} \left|{\frac{\partial \varphi}{\partial x_2}(x_1, x_2, x_3)}\right|^2\,d x_1\,d x_2\,d x_3
\\\nonumber= \beta_f \int_{-s_0}^{s_0}\int_{{\{x_1^2+ x_2^2\le s_0^2- x_3^2\}}}\left|{i \frac{\partial \varphi}{\partial x_1}(x_1, x_2, x_3)+A_1(x_1, x_2, x_3) \varphi(x_1, x_2, x_3)}\right|^2\,d x_1\,d x_2\,d x_3\\\nonumber
+\beta_f \int_{-s_0}^{s_0}\int_{{\{x_1^2+x_2^2\le s_0^2-x_3^2\}}} \left|{\frac{\partial \varphi}{\partial x_2}(x_1, x_2, x_3)}\right|^2\,d x_1\,d x_2\,d x_3\\\nonumber
\ge\beta_f\int_{-s_0/\sqrt{2}}^{s_0/\sqrt{2}}\int_{{\{x_1^2+ x_2^2\le s_0^2- x_3^2\}}}\left|{i \frac{\partial \varphi}{\partial x_1}(x_1, x_2, x_3)+ A_1(x_1, x_2, x_3) \varphi(x_1, x_2, x_3)}\right|^2\,d x_1\,d x_2\,d x_3\\\nonumber
+\beta_f \int_{-s_0/\sqrt{2}}^{s_0/\sqrt{2}}\int_{{\{x_1^2+ x_2^2\le s_0^2- x_3^2\}}} \left|{\frac{\partial \varphi}{\partial x_2}(x_1, x_2, x_3)}\right|^2\,d x_1\,d x_2\,d x_3\\\nonumber
=\beta_f \int_{-s_0/\sqrt{2}}^{s_0/\sqrt{2}}\int_{{\{x_1^2+ x_2^2\le s_0^2- x_3^2\}}}\biggl(\left|{i \frac{\partial \varphi}{\partial x_1}(x_1, x_2, x_3)+ A_1(x_1, x_2, x_3) \varphi(x_1, x_2, x_3)}\right|^2
\\\label{asympt.}+\left|{\frac{\partial \varphi}{\partial x_2}(x_1, x_2, x_3)}\right|^2\biggr)\,d x_1\,d x_2\,d x_3\,.
\end{eqnarray}

We will use the following Lemma to continue the proof. The proof of this Lemma is given in the Appendix.  

\begin{lemma}\label{Lemma}
Let us denote by $\lambda_1(\tilde{B}, R)$ the ground state eigenvalue of the Neumann magnetic Laplacian with a constant magnetic field 
$\tilde{B}$ defined on a ball with center at the origin and radius $R>\epsilon,\,\epsilon>0$. Then, for large values of $\tilde{B}$, uniformly with respect to $R$ the following asymptotic takes place
$$
\lambda_1(\tilde{B}, R)= \alpha \tilde{B}+ o(\tilde{B}) 
$$
where $\alpha$ is some positive constant independent on $\tilde{B}$ and $R$.
\end{lemma}
From the above Lemma, it follows that there exists some $\gamma>0$ such that if $\tilde{B}\ge \gamma$, then for any $R\in \left(\frac{s_0}{\sqrt{2}}, s_0\right)$ one has 
\begin{equation}\label{cond.}
\lambda_1\left(\tilde{B}, R\right)\ge \frac{\alpha \tilde{B}}{2}\,.
\end{equation}
In view of the fact that the divergence of the magnetic field $B$ must to be equal to zero, 
we can conclude that for any fixed $x_3\in \left(-\frac{s_0}{\sqrt{2}}, \frac{s_0}{\sqrt{2}}\right)$, the two-dimensional magnetic field corresponding to the magnetic potential $(A_1(x_1, x_2, x_3), 0)$ defined on the two-dimensional ball $\{x_1^2+ x_2^2\le s_0^2- x_3^2\}$ equals to $B_3(x_1, x_2, x_3)=B_3^0$. 

Hence
\begin{eqnarray*}
\int_{{\{x_1^2+ x_2^2\le s_0^2- x_3^2\}}}\left(\left|{\frac{\partial \varphi}{\partial x_1}(x_1, x_2, x_3)+ A_1(x_1, x_2, x_3) \varphi(x_1, x_2, x_3) }\right|^2
+\left|{\frac{\partial \varphi}{\partial x_2}(x_1, x_2, x_3)}\right|^2\right)\,d x_2\,d x_3
\\ \ge\frac{\alpha { B_3^0}}{2} \int_{{\{x_1^2+ x_2^2\le s_0^2- x_3^2\}}}|\varphi(x_1, x_2, x_3)|^2\,d x_2\,d x_3\,.
\end{eqnarray*}
Applying the above inequality to the right-hand side of the estimate (\ref{asympt.}), one obtains
\begin{eqnarray*}
\int_{\mathbb{R}^3}\left|i \nabla \varphi(x_1, x_2, x_3) + A(x_1, x_2, x_3) \varphi(x_1, x_2, x_3)\right|^2 f(x_1, x_2)^2\,d x_1\,d x_2\,d x_3\\\ge
\frac{\alpha \beta_f {B_3^0}}{2}  \int_{-s_0/\sqrt{2}}^{s_0/\sqrt{2}}\int_{{\{x_1^2+ x_2^2\le s_0^2- x_3^2\}}} |\varphi(x_1, x_2, x_3)|^2\,d x_1\,d x_2\,d x_3\,.
\end{eqnarray*}

Returning to inequality (\ref{est.}) and using the above bound, we obtain that for values of ${B_3^0}$ starting from $\gamma>0$
\begin{eqnarray}\nonumber
\int_{\mathbb{R}^3}\left(|i \nabla \psi(x_1, x_2, x_3)+ A(x_1, x_2, x_3) \psi(x_1, x_2, x_3)|^2- V(x_1, x_2, x_3) |\psi(x_1, x_2, x_3)|^2\right)\,d x_1\,d x_2\,d x_3 
\\\nonumber\ge \frac{\alpha \beta_f {B_3^0}}{2} \int_{-s_0/\sqrt{2}}^{s_0/\sqrt{2}}\int_{{\{x_1^2+ x_2^2\le s_0^2- x_3^2\}}} |\varphi(x_1, x_2, x_3)|^2\,d x_1\,d x_2\,d x_3\\\nonumber- \int_{\mathbb{R}^3}(\tilde{V}(x_1, x_2, x_3)-V(x_1, x_2)) f(x_1, x_2)^2 |\varphi(x_1, x_2, x_3)|^2\,d x_1\,d x_2\,d x_3\\\nonumber+ e \int_{\mathbb{R}^3} |\psi(x_1, x_2, x_3)|^2\,d x_1\,d x_2\,d x_3\\\nonumber= \frac{\alpha \beta_f {B_3^0}}{2} \int_{\mathcal{B}(0, s_0)\cap\left\{{|x_3|}\le\frac{s_0}{\sqrt{2}}\right\}} |\varphi(x_1, x_2, x_3)|^2\,d x_1\,d x_2\,d x_3\\\nonumber- \int_{\mathbb{R}^3}(\tilde{V}(x_1, x_2, x_3)-V(x_1, x_2)) f(x_1, x_2)^2 |\varphi(x_1, x_2, x_3)|^2\,d x_1\,d x_2\,d x_3\\\nonumber+ e \int_{\mathbb{R}^3} |\psi(x_1, x_2, x_3)|^2\,d x_1\,d x_2\,d x_3\\\nonumber=\frac{\alpha \beta_f {B_3^0}}{2} \int_{\mathcal{B}(0, s_0)\cap\left\{{|x_3|}\le\frac{s_0}{\sqrt{2}}\right\}} |\varphi(x_1, x_2, x_3)|^2\,d x_1\,d x_2\,d x_3\\\nonumber- \int_{\mathcal{B}(0, s_0)
\cap\left\{|{x_3}|\le\frac{s_0}{\sqrt{2}}\right\}}(\tilde{V}(x_1, x_2, x_3)-V(x_1, x_2)) f(x_1, x_2)^2 |\varphi(x_1, x_2, x_3)|^2\,d x_1\,d x_2\,d x_3\\\nonumber-
 \int_{\mathbb{R}^3\setminus\left(\mathcal{B}(0, s_0)\cap\left\{{|x_3|}\le\frac{s_0}{\sqrt{2}}\right\}\right)}(\tilde{V}(x_1, x_2, x_3)-V(x_1, x_2)) f(x_1, x_2)^2 |\varphi(x_1, x_2, x_3)|^2\,d x_1\,d x_2\,d x_3\\\nonumber+ e \int_{\mathbb{R}^3} |\psi(x_1, x_2, x_3)|^2\,d x_1\,d x_2\,d x_3\,.
\end{eqnarray}

In view of the construction of $\tilde{V}$ and using the non-negativeness of $V$
the above inequality takes the following form
\begin{eqnarray}\nonumber
\int_{\mathbb{R}^3}\left(|i \nabla \psi(x_1, x_2, x_3)+ A(x_1, x_2, x_3) \psi(x_1, x_2, x_3)|^2- V(x_1, x_2, x_3) |\psi(x_1, x_2, x_3)|^2\right)\,d x_1\,d x_2\,d x_3\\\nonumber\ge\frac{\alpha  \beta_f B_3^0}{2} \int_{\mathcal{B}(0, s_0)\cap\left\{|x_3|\le\frac{s_0}{\sqrt{2}}\right\}} |\varphi(x_1, x_2, x_3)|^2\,d x_1\,d x_2\,d x_3\\\nonumber- \int_{\mathcal{B}(0, s_0)\cap\left\{|x_3|\le\frac{s_0}{\sqrt{2}}\right\}}(\tilde{V}(x_1, x_2, x_3)-V(x_1, x_2)) f(x_1, x_2)^2 |\varphi(x_1, x_2, x_3)|^2\,d x_1\,d x_2\,d x_3\\\label{fin.}+ e \int_{\mathbb{R}^3} |\psi(x_1, x_2, x_3)|^2\,d x_1\,d x_2\,d x_3\,.
\end{eqnarray}

Therefore, by choosing this assumption
$$
\|\tilde{V}(x_1, x_2, x_3)-V(x_1, x_2)\|_{L^\infty\left(\mathcal{B}(0, s_0)\cap\left\{|{x_3}|\le s_0/\sqrt{2}\right\}\right)} \|f\|_{L^\infty(\mathbb{R}^2)}^2
\le \frac{\alpha \beta_f {B_3^0}}{2}
$$
the right-hand side of inequality (\ref{fin.}) will be no less than $e$. It is easy to verify that this assumption is guaranteed by (\ref{assumption2}) with the constant
$$
C= \frac{\alpha \beta_f}{2 \|f\|_{L^\infty(\mathbb{R}^2)}^2}.
$$

\end{proof}

\section{Appendix}
\setcounter{equation}{0}

In this section, we will prove the Lemma \ref{Lemma}. Let $\mathcal{B}(0, R)$ be the two-dimensional disk with centered at the origin with radius $R>\epsilon$ and let $\tilde{A}$ be the magnetic potential corresponding to the constant magnetic field $\tilde{B}$. After making a change of variables in the quadratic form of the magnetic Laplacian with the magnetic potential $\tilde{A}$, defined on the $\mathcal{B}(0, R)$, we have
$$\int_{\mathcal{B}(0, R)}|i \nabla u+ \tilde{A} u|^2\,d z_1\,d z_2= \int_{\mathcal{B}(0, 1)}|i \nabla v+ \tilde{A}^1 v|^2\,d x_1\,d x_2\,,$$
where $\mathcal{B}(0, 1)$ is the unit two-dimensional disk, $\tilde{A}^1= R \tilde{A}(R x_1, Rx_2)$, and $v= u(R x_1, R x_2)$. 

Since the magnetic field corresponding to the magnetic potential $R \tilde{A}(R x_1, R x_2)$ is constant and equal to $R^2 \tilde{B}$ then in view of the works \cite{S19} and later in \cite{KM25}, \cite{HL25}, the ground state eigenvalue $\lambda_1(R^2 \tilde{B}, 1)$ of the magnetic Neumann Laplacian on a unit two-dimensional disk and magnetic field $R^2 \tilde{B}$  has the following asymptotic behaviour
$$
\lambda_1(R^2 \tilde{B}, 1)= \alpha R^2 \tilde{B}+ o(R^2 \tilde{B}) 
$$
where $\alpha$ is a positive constant that is independent of both $\tilde{B}$ and $R$. Note that the above asymptotic behaviour is uniform with respect to $R$ in view of $R>\epsilon$. Since $\int_{\mathcal{B}(0, 1)}|v|^2\,d x_1\,d x_2=\frac{1}{R^2} \int_{\mathcal{B}(0, R)}|u|^2\,d z_1\,d z_2$,
this proves the Lemma.

\begin{remark}\label{final}

In the paper \cite{E22}, it was proven that if $\Gamma$ is a $C^2$-smooth curve, which is a local deformation of the straight line $\{(0, 0, x_3)\}_{x_3\in\mathbb{R}}$ satisfying the so-called Tang condition, then the operator $-\Delta- \tilde{V}$ with $\tilde{V}$ given by (\ref{potential}) and  $V= \frac{1}{\varepsilon}\chi_\omega$, where $\omega$ is some open subset in $\mathbb{R}^2$, has a non-empty discrete spectrum provided that $\varepsilon$ is small enough.

Assume that assumption (\ref{assumption1}) is satisfied. Our results show that by adding the suitable magnetic field which satisfies assumption (\ref{assumption2}) simply by choosing 
${B_3^0}\ge \frac{2}{C \varepsilon}$ the discrete spectrum of the operator $(i \nabla+A)^2- \tilde{V}$ becomes empty.

\end{remark}

\section*{Declarations}

\subsection*{Funding} This work is supported in part by Instituto Polit\'ecnico Nacional (grant number SIP20253412).

\subsection*{Competing Interests} The authors declare that they have no competing interests regarding the publication of this paper.

\subsection*{Author contributions} All authors contributed equally to the study, read and approved the final version of the submitted manuscript.

\subsection*{Availability of data and material} Not applicable

\subsection*{Code availability} Not applicable

\end{document}